\documentclass[11pt,a4paper,reqno,centertags]{article}
\usepackage{xypic,a4wide,amsmath,amsthm,amssymb}

\newtheorem{thm}{Theorem}[section]

\theoremstyle{remark}
\newtheorem{rem}[thm]{Remark}
\newtheorem{ex}[thm]{Example}

\newenvironment{example}{\begin{ex}\rm}{\qee\end{ex}}
\newcommand{\qee}{\mbox{\hspace{0.2mm}}\hfill$\triangle$}

{\nonumber}
\newcommand{\di}{{\mathrm d}}
\newcommand{\R}{{\mathbb R}}
\newcommand{\C}{{\mathbb C}}

\begin{document}
\begin{center}
 {\LARGE\bf Gromov-Witten invariants of target curves via Symplectic Field Theory.}\\[15pt]
 {\sc Paolo Rossi}\\
 {(SISSA - Trieste)}
\end{center}

\vspace{1cm}

\begin{abstract}
We compute the Gromov-Witten potential at all genera of target smooth Riemann surfaces using Symplectic Field Theory techniques and establish differential equations for the full descendant potential. 
This amounts to impose (and possibly solve) different kinds of Schr\"odinger equations related to some quantization of the dispersionless KdV hierarchy. In particular we find very explicit formulas for the Gromov-Witten invariants of low degree of $\mathbb{P}^1$ with descendants of the K\"ahler class.
\end{abstract}

\section*{Introduction}

Symplectic Field Theory \cite{EGH} is a fairly new branch of symplectic topology studying holomorphic curves in symplectic manifolds with ends in the spirit of Gromov-Witten theory (which it actually contains as a special case). The presence of cylindrical ends that are symplectizations of contact manifolds and, in particular, the imposition (boundary condition) that holomorphic curves have punctures that asymptotically coincide with Reeb orbits in these ends, gives to the theory a rich algebraic structure where classical and quantum integrable systems, together with many natural tools of that context, arise.

In particular let $W$ be a symplectic manifold with cylindrical ends, i.e. a symplectic cobordism between contact manifolds $V^+$, $V^-$, completed by attaching to these boundaries their symplectizations $V^+\times [0,+\infty)$, $V^-\times (-\infty,0]$. The potential counting holomorphic curves in each of the two symplectizations $V^+\times \R$, $V^-\times \R$, asymptotically cylindrical over Reeb orbits of $V^+$, $V^-$, is to be interpreted as the Hamiltonian for a (quantum) Hamiltonian system. Schr\"odinger equation (and its semiclassical limit, Hamilton-Jacobi equation) relative to this Hamiltonian enters in the computation of the Gromov-Witten potential of $W$, that in this context plays the role of the phase of the wave function, as explained in \cite{EGH}, section $2.7$. Moreover analogous formulas (see \cite{EGH}, section $2.5$) can be used to deduce potentials for a composition of cobordisms from the ones of the single pieces.

Symplectic topologists are fairly confident that the Hamiltonian systems arising this way in Symplectic Field Theory are in general integrable systems (it is so for all explicitly known cases, which are indeed quite few). Even more, it is believed that the other Hamiltonians of the integrable hierarchy, and their correponding Schr\"odinger equations, can be used in a way totally analogous to the one explained above to compute the descendant potential of $W$.

\vspace{1cm}

In this paper we use these Symplectic Field Theory techniques to compute at all genera the Gromov-Witten potential of target curves of any genus. Our program is similar to the one via relative Gromov-Witten invariants of Okounkov and Pandharipande (\cite{OP2}), i.e. we first consider the two building blocks of a Riemann surface (the cap and the pair of pants) and then attach them to obtain the Gromov-Witten potential of the curve $\Sigma_g$ of genus $g$. Indeed Symplectic Field Theory can be interpreted as a theory of relative invariants, so the analogy in this sense is strong. Nonetheless our method of computation is not based on Virasoro constraints (as for \cite{OP2}), but on relations proved in \cite{EGH} for a very general setting, which makes the computation not only elegant, but even fairly explicit. Moreover, thanks to the surprising emergence of integrable systems from the formalism of Symplectic Field Theory, we can give a beautiful interpretation of many aspects of the full descendants GW-potential in terms of a quantization of the dispersionless KdV hierarchy (after an idea of Eliashberg, \cite{E}).

\subsection*{Acknowledgements} 
I am very grateful to Y. Eliashberg and A. Givental for suggestions and a fruitful discussion (particularly regarding formula (\ref{qdkdvgenfun})). Moreover I wish to thank my advisor B. Dubrovin for his constant support and expert guidance.

\vspace{0.5cm}

This work is partially supported by the European Science Foundation Programme ``Methods of Integrable Systems, Geometry, Applied Mathematics" (MISGAM), the Marie Curie RTN ``European Network in Geometry, Mathematical Physics and Applications"  (ENIGMA),  and by the Italian Ministry of Universities and Researches (MIUR) research grant PRIN 2004 ``Geometric methods in the theory of nonlinear waves and their applications".

\section{Symplectic Field Theory preliminaries}

For the reader's convenience and also in order to fix notations, we recall here from \cite{EGH} the main Symplectic Field Theory tools hinted above. Our technique will basically consist in constructing the Riemann surface by attaching simpler symplectic cobordisms, namely the complex line $\C$ (or the cap) and $\mathbb{P}^1-\{0,1,\infty\}$ (or the pair of pants). This kind of computations will require also to consider the cylinder $S^1 \times \R$, as already done in \cite{EGH}. 

First recall that, to the symplectization $(V\times \R, \di (\mathrm{e}^t \alpha))$ of a contact manifold $(V,\alpha)$, we associate the Hamiltonian (SFT-potential) $$\mathbf{H}=\frac{1}{\hbar}\sum_{g=0}^{\infty}\mathbf{H}_g \hbar^g$$ where $$\mathbf{H}_g=\sum_{A\in H_2(V)} \ \sum_{r,s^\pm=0}^\infty \frac{1}{r!s^+ !s^- !} \langle t,\ldots,t;q,\ldots,q;p,\ldots,p \rangle^A_g \ z^A$$
Here $t=\sum t_i \Theta_i \in \Omega^\bullet(V)$, while $p=\sum \frac{1}{\mu_\gamma} p_\gamma [\gamma]$, $q=\sum \frac{1}{\mu_\gamma} q_\gamma [\gamma]$ (the sum is over the set $\mathcal{P}(V)$ of periodic Reeb orbits, including multiples with multiplicity $\mu_\gamma$) and
\begin{equation*}
\begin{split}
&\langle \Theta_1,\ldots,\Theta_r;\gamma^-_1,\ldots,\gamma^-_{s^-};\gamma^+_1,\ldots,\gamma^+_{s^+} \rangle^A_g =\\
&\int_{\mathcal{M}^A_{g,r,s^-,s^+}/\R} \mathrm{ev}^*(\Theta_1\otimes\ldots\otimes\Theta_r\otimes\gamma^-_1\otimes\ldots\otimes\gamma^-_{s^-}\otimes \gamma^+_1\otimes\ldots\otimes\gamma^+_{s^+}) 
\end{split}
\end{equation*}
The integral is over the moduli space of holomorphic curves in $V\times \R$ with $r$ marked points and $s^\pm$ positive/negative punctures asymptotically cylindrical over Reeb orbits and realizing, together with the chosen capping surfaces (see \cite{EGH} for details), the homology cycle $A$ in $V\times \R$, modulo the $\R$ action coming from the $\R$ symmetry of the cylindrical target space $V\times \R$. The map $\mathrm{ev}$ is the natural evaluation map $\mathrm{ev}:\mathcal{M}^A_{g,r,s^-,s^+}\to V^r\times\mathcal{P}(V)^{s^-}\times\mathcal{P}(V)^{s^+}$.

Similarly, to a general completed symplectic cobordism $W=\overrightarrow{V^+ V^-}$, we associate the SFT-potential $$\mathbf{F}=\frac{1}{\hbar}\sum_{g=0}^{\infty}\mathbf{F}_g \hbar^g$$ where $$\mathbf{F}_g=\sum_{A\in H_2(W)} \ \sum_{r,s^\pm=0}^\infty \frac{1}{r!s^+ !s^- !} \langle t,\ldots,t;q,\ldots,q;p,\ldots,p \rangle^A_g \ z^A$$
Here $t=\sum t_i \Theta_i \in \Omega^\bullet(W)$, while $p=\sum \frac{1}{\mu_\gamma} p_\gamma [\gamma]$, $q=\sum \frac{1}{\mu_\gamma} q_\gamma [\gamma]$ (the sums are respectively over the sets $\mathcal{P}(V^+)$ and $\mathcal{P}(V^-)$ of periodic Reeb orbits, including multiples with multiplicity $\mu_\gamma$) and
\begin{equation*}
\begin{split}
&\langle \Theta_1,\ldots,\Theta_r;\gamma^-_1,\ldots,\gamma^-_{s^-};\gamma^+_1,\ldots,\gamma^+_{s^+} \rangle^A_g =\\
&\int_{\mathcal{M}^A_{g,r,s^-,s^+}} \mathrm{ev}^*(\Theta_1\otimes\ldots\otimes\Theta_r\otimes\gamma^-_1\otimes\ldots\otimes\gamma^-_{s^-}\otimes \gamma^+_1\otimes\ldots\otimes\gamma^+_{s^+}) 
\end{split}
\end{equation*}
The integral is over the moduli space of holomorphic curves in $W$ with $r$ marked points and $s^\pm$ positive/negative punctures asymptotically cylindrical over Reeb orbits and realizing, together withe the chosen capping surfaces, the homology cycle $A$ in $W$. The map $\mathrm{ev}$ is the natural evaluation map $\mathrm{ev}:\mathcal{M}^A_{g,r,s^-,s^+}\to W^r\times\mathcal{P}(V^-)^{s^-}\times\mathcal{P}(V^+)^{s^+}$.

If we assign the following grading to the variables ($2n$ is the dimension of the symplectic cobordism and $\mathrm{CZ}(\gamma)$ is the Conley-Zehnder index of $\gamma$):
\begin{align*}
&\mathrm{deg}(t_i)=\mathrm{deg}(\Theta_i)-2 \hspace{1cm} &\mathrm{deg}(p_\gamma)=-\mathrm{CZ}(\gamma)+(n-3) \\
&\mathrm{deg}(\hbar)=2(n-3) & \mathrm{deg}(q_\gamma)=+\mathrm{CZ}(\gamma)+(n-3) \\
&\mathrm{deg}(z^A)=-2c_1(A)
\end{align*}
(where $c_1(A)$ is the first Chern class of $\text{T}W|_A$) then the potentials $\mathbf{H}$ and $\mathbf{F}$ can be see as elements in graded Weyl algebras where all the symbols supercommute except $p_\gamma$ with $q_\gamma$ (same periodic orbit $\gamma$) for which the supercommutator is $[p_\gamma,q_\gamma]=\mu_\gamma \hbar$. Inside this graded Weyl algebra the function $\mathbf{H}$ satisfies a structure equation in the form $\di\mathbf{H}+\frac{1}{2}[\mathbf{H},\mathbf{H}]=0$ (see \cite{EGH}).

\vspace{1cm}

In computing these potentials, it is essential to know the (virtual) dimension of the moduli space $\mathcal{M}^A_{g,r,s^-,s^+}$, that is given in \cite{EGH} by the index formula $$\text{dim} \mathcal{M}^A_{g,r,s^-,s^+}=\sum_{1}^{s^+}\text{CZ}(\gamma_i^+)+\sum_1^{s^-}\text{CZ}(\gamma_k^-)+(n-3)(2-2g-s^+ -s^-)+2c_1^\text{rel}(A_\text{rel})+2r$$
where this time $c_1^\text{rel}(A_\text{rel})$ is the first relative Chern class of $\text{T}W|_{A_\text{rel}}$ and $A_\text{rel}$ is the relative homology cycle realized by the holomorphic curves in $W$.
\vspace{1cm}

The first result from \cite{EGH} we need to use is the following. Let us assume that $W$ has only a positive end $V\times[0,+\infty)$, and choose what is called a \emph{basic system} $\Delta_1,\dots,\Delta_k$, $\Theta_1,\dots,\Theta_m$ of closed differential forms on $W$, with cylindrical ends, such that
\begin{itemize}
\item{a)} $\Delta_1,\dots,\Delta_k \in H^\bullet(W)$, and the restrictions $\delta_i=\Delta_i|_V$, $i=1,\dots,l$ for $l\leq k$ are independent elements in $\mathrm{Im}(H^\bullet(W)\to H^\bullet(V))$;
\item{b)}  $\Theta_1,\dots,\Theta_m$ are compactly supported and independent elements in $\mathrm{Ker}(H^\bullet_{\mathrm{comp}}(W)\to H^\bullet(W))$,
\item{c)} there exist forms $\theta_1,\dots,\theta_m$ on $V$ and a compactly supported $1$-form $\rho$ on $(0,+\infty)$, such that $\Theta_j=\rho\wedge \theta_j,\;j=1,\dots, m$.
\end{itemize}

\begin{thm}[\cite{EGH}]\label{Schrodinger}
Let $\mathbf{H}$ be the Hamiltonian associated with the
contact manifold $V$.
Set $$\mathbf{H}^j(t_1,\dots,t_l,q,p)=\Big(\frac{\partial \mathbf{H}}{\partial s_j}(\sum\limits_{i=1}^l t_i\delta_i+s_j\theta_j,q,p)\Big)\Big|_{s_j=0},\;j=1,\dots m,$$
$$\mathbf{F}^0(t_1,\dots,t_k, p)=\mathbf{F}(\sum t_i\Delta_i,p).$$
Then the potential associated to $W$ is given by:
$$\mathrm{e}^{\mathbf{F}(\sum t_i\Delta_i+\sum \tau_j \Theta_j,p)}=\mathrm{e}^{\mathbf{F}^0(t_1,\dots,t_k,p)} \mathrm{e}^{\tau_m\overleftarrow{\mathbf{H}^m}(t_1,\dots,t_l,q,p)}\dots \mathrm{e}^{\tau_1\overleftarrow{\mathbf{H}^1}(t_1,\dots,t_l,q,p)},$$ where $\overleftarrow{\mathbf{H}^i}$ is the operator obtained from $\mathbf{H}^i$ by quantizing $q_\gamma=\mu_\gamma \hbar\overleftarrow{\frac{\partial }{\partial p_\gamma}}$.
\end{thm}

\vspace{1cm}

In what follows we also make use of the following theorem, concerning the composition $W=\overrightarrow{V^-V^+}$ of two symplectic cobordisms $W_-=\overrightarrow{V^- V}$ and $W_+=\overrightarrow{VV^+}$.

\begin{thm}[\cite{EGH}]\label{composition} 
Let us denote by $\mathbf{F}_W$, $\mathbf{F}_{W_-}$ and $\mathbf{F}_{W_+}$ the SFT-potentials of $W$, $W_-$ and $W_+$ respectively. Notice that any cohomology class in $H^\bullet(W)$ can be represented by a form $t$ which splits into the sum of forms $t_\pm$ with cylindrical ends on $W_\pm$ so that $t_\pm|_V=t_V$. Then:
$$\mathbf{F}_W(q^-,p^+,t)=\mathbf{F}_{W_-}(q^-,p,t_-)\Diamond \mathbf{F}_{W_+}(q,p^+,t_+)$$
where $\displaystyle{\mathrm{e}^{F \Diamond G}=\left.\left(\mathrm{e}^{\overrightarrow{F}} \mathrm{e}^{G}\right)\right|_{q=0}}$ and $\overrightarrow{F}$ is the operator obtained by quantizing $\displaystyle{p_{\gamma}=\hbar \mu_{\gamma} z^{A_\gamma} \overrightarrow{\frac{\partial}{\partial q_{\gamma}}}}$ and $A_\gamma\in H_2(W)$ is the cycle in $W$ formed by the capping surfaces (see again \cite{EGH} for details) of $\gamma$ in $W_-$ and $W_+$.
\end{thm}

\section{Pair of pants potential}

In this section we compute the SFT-potential $\mathbf{F}_\mathrm{pants}$ for the \emph{pair of pants}, i.e. the completed symplectic cobordism $W$ between $V^-=S^1$ and $V^+=S^1 \coprod S^1$. $W$ can also be seen as the complex projective line $\mathbb{P}^1$ minus $\{0,1,\infty\}$, with its standard K\"ahler structure. This way $V^+$ corresponds to the circles around, say, $0$ and $\infty$, and $V^-$ to the circle around $1$. Call then $\phi$ the longitude on the Riemann sphere and $\varphi$ its restriction to $V^+$; call $\alpha$ the angle on the circle $V^-$ around $1$.

Choose the basic system of forms as
$$\Delta_0=1 \in H^0(W)\quad \text{restricting to}\quad \delta_0^+=(1,1) \in H^0(V^+),\quad \delta_0^-=1 \in H^0(V^-),$$
$$\Delta_1=\di \phi \in H^1(W)\quad \text{restricting to}\quad \delta_1^+=(\di \varphi,\di \varphi) \in H^1(V^+),\quad \delta_1^-=0 \in H^1(V^-)$$
and
$$\Theta_1 \in H^1_{\mathrm{comp}}(W)\quad \text{projecting to}\quad \theta_1^+=(1,-1) \in H^0(V^+),\quad \theta_1^-=0 \in H^0(V^-)$$
$$\Theta_2 \in H^2_{\mathrm{comp}}(W)\quad \text{projecting to}\quad \theta_2^+=(0,0) \in H^1(V^+),\quad \theta_2^-=\di \alpha \in H^1(V^-).$$

\vspace{1cm}

The SFT-potential of $V^+=S^1 \coprod S^1$ (or its symplectization) at all genera is then easily computed from the one of $S^1 \times \R$ (see \cite{EGH}). What we will need is:
$$\mathbf{H}^+(\sum t_i \delta_i^+ + \sum s_j \theta_j^+)=\frac{1}{\hbar}\left[t_0^2t_1 +t_0t_1s_1+s_1^2t_1 + t_1\left(\sum q^1_k p^1_k + q^2_k p^2_k \right) -\frac{\hbar t_1}{12} \right]$$
and
\begin{equation}\label{cylinder}
\mathbf{H}^-(\sum t_i \delta_i^- + \sum s_j \theta_j^-)=\frac{1}{\hbar}\left[\frac{t_0^2s_2}{2} + s_2\sum q_k p_k-\frac{\hbar s_2}{24} \right]
\end{equation}
where $\mathrm{deg}t_1=\mathrm{deg}s_2=-1$, $\mathrm{deg}t_0=\mathrm{deg}s_1=-2$, $\mathrm{deg}p^i_k=\mathrm{deg}p_k=\mathrm{deg}q^i_k=\mathrm{deg}q_k=-2$ and $\mathrm{deg}\hbar=-4$.

Define
\begin{align}
\mathbf{H}^+_1:=& \left.\frac{\partial \mathbf{H^+}}{\partial s_1}\right|_{s_i=0}=\frac{1}{\hbar}t_0t_1 \\
\mathbf{H}^-_2:=& \left.\frac{\partial \mathbf{H^-}}{\partial s_2}\right|_{s_i=0}=\frac{1}{\hbar} \left[\frac{t_0^2}{2} + \sum q_k p_k - \frac{\hbar}{24}\right]
\end{align}
We will refer to these functions as Hamiltonians. Notice, by the way, that putting 
\begin{equation}\label{fourier}
u(x)=t_0+ \sum p_k \mathrm{e}^{\mathrm{i}kx} + q_k \mathrm{e}^{-\mathrm{i}kx}
\end{equation}
we have $$\mathbf{H}^-_2=\frac{1}{\hbar}:\int_{S^1}\left(\frac{u(x)^2}{2} -\frac{\hbar}{24}\right)\di x:$$
where the normal ordering $:\cdot:$ means that the $q$ variables are to be put on the left of the $p$ variables as if, inside the colon symbols, $q$'s and $p$'s all commuted.

We use these Hamiltonians and Theorem \ref{Schrodinger} on the initial datum $\mathbf{F}_\text{pants}(\sum t_i \Delta_i)$ which is computed directly, by dimension counting. In fact from a combination of the reconstruction theorem for ramified coverings of $\mathbb{P}^1$, the Riemann-Hurwitz theorem and the index formula for the dimension of $\mathcal{M}^A_{g,r,s^+,s^-}$, we get the known relation (see e.g. \cite{OP1}) between relative Gromov-Witten (or SFT-) potential and Hurwitz numbers
\begin{equation}
\mathbf{F}_\text{pants}(\sum t_i \Delta_i)=\sum_{g}\sum_{d}\sum_{|\mu^0|,|\mu^1|,|\mu^\infty|=d} H^{\mathbb{P}^1}_{g,d}(\mu^0,\mu^1,\mu^\infty)(p^1)^{\mu^0}(p^2)^{\mu^\infty}q^{\mu^1}\hbar^{g-1}
\end{equation}
where $H^{\mathbb{P}^1}_{g,d}(\mu^0,\mu^1,\mu^\infty)$ is the Hurwitz number counting coverings of $\mathbb{P}^1$ of degree $d$ and genus $g$, branched only over $0$, $1$, and $\infty$ with ramification profile given by $\mu^0$, $\mu^1$ and $\mu^\infty$ respectively.

Then, by Theorem \ref{Schrodinger} and recalling that $\mathrm{e}^{sx\partial_x}f(x)=f(x\mathrm{e}^s)$, one gets
\begin{equation}\label{pairofpants}
\begin{split}
\mathbf{F}_\text{pants}&(\sum t_i \Delta_i+\sum s_j \Theta_j)= \text{log}\left(\mathrm{e}^{s_2 \overrightarrow{\mathbf{H}^-_2}}\, \mathrm{e}^{\mathbf{F}_\text{pants}(t)} \, \mathrm{e}^{s_1 \mathbf{H}^+_1}\right)=
\frac{1}{\hbar}\left[\frac{t_0^2s_2}{2}+t_0t_1s_1
-\frac{\hbar s_2}{24}\right]\\
&+\sum_{g}\sum_{d}\sum_{|\mu^0|,|\mu^1|,|\mu^\infty|=d} H^{\mathbb{P}^1}_{g,d}(\mu^0,\mu^1,\mu^\infty)(p^1)^{\mu^0}(p^2)^{\mu^\infty}\left(q\mathrm{e}^{s_2}\right)^{\mu^1}\hbar^{g-1}
\end{split}
\end{equation}
where $\displaystyle{\overrightarrow{\mathbf{H}^-_2}}$ means the operator obtained from $\mathbf{H}^-_2$ by quantizing $\displaystyle{p_k=k \hbar \overrightarrow{\frac{\partial}{\partial q_k}}}$.

\vspace{1cm}

From this potential, together with the one for the cap (see e.g. \cite{EGH})
\begin{equation}\label{cap}
\mathbf{F}_\text{cap}(t_0 1+s_2 (\rho \wedge \di \phi))= \frac{1}{\hbar}\left[\frac{t_0^2s_2}{2}-\frac{s_2}{24}+\frac{1}{\hbar}\mathrm{e}^{s_2}p_1\right]
\end{equation}
and using Theorem \ref{composition} one gets the Gromov-Witten potential for the general genus $g$ Riemann surface. In the next section we illustrate this methods of computation for the simplest case of the elliptic curve $E$, for which the result is very explicit.

\section{Gromov Witten potential of $E$}

Let $E$ be a smooth elliptic curve. Consider its pair-of-pants decomposition; it consists of two caps and two pairs of pants. Let $U$ be one of the halves of the torus $E$ when cut along two representatives of the same element in $\pi_1(E)$. $U$ is the symplectic manifold with a contact boundary resulting from attaching a cap to a pair of pants. Its potential is then given by (\ref{pairofpants}) and (\ref{cap}) by applying Theorem \ref{composition}:
\begin{equation*}
\begin{split}
\mathbf{F}_U(\sum t_i \Delta_i+\sum s_j \Theta_j)&=\mathrm{log}\left.\left(\mathrm{e}^{\overrightarrow{\mathbf{F}_\text{cap}}(0)}\quad \mathrm{e}^{\mathbf{F}_\text{pants}(\sum t_i \Delta_i+\sum s_j \Theta_j)}\right)\right|_{q=0}\\
&=\mathrm{log}\left.\left(\mathrm{e}^{\overrightarrow{\partial_{q_1}}}\quad \mathrm{e}^{\mathbf{F}_\text{pants}(t,s)}\right)\right|_{q=0}
\end{split}
\end{equation*}
This results in a substantial simplification of the term containing Hurwitz numbers in (\ref{pairofpants}) since putting $q_i=0$ after derivation with respect to $q_1$ selects just the terms counting for coverings branched only over $0$ and $+\infty$, hence with the same branching number (as it is natural for the potential of the space $U$). So we get
\begin{equation}
\mathbf{F}_U(\sum t_i \Delta_i+\sum s_j \Theta_j)=\frac{1}{\hbar}\left[\frac{t_0^2s_2}{2}+t_0t_1s_1-\frac{\hbar s_2}{24}+\sum\frac{1}{k}p^1_kp^2_k\mathrm{e}^{ks}\right]
\end{equation}
The last step towards the computation of the Gromov-Witten potential of $E$ is performed by gluing two copies of $U$, one with positive ends and one with negative ends. Let $\tilde{\Delta}_i, \tilde{\Theta}_j \in H^\bullet(E)$ be the obvious extensions of $\Delta_i, \Theta_j$ from $U$ to the entire $E$, then the gluing Theorem \ref{composition} takes the form:
$$\mathcal{F}_E \left(\sum t_i \tilde{\Delta}_i+ \sum s_j \tilde{\Theta}_j\right)=\left.\left(\mathrm{e}^{\overrightarrow{\mathbf{F}_{U^+}}(t,s)}\quad \mathrm{e}^{\mathbf{F}_{U^-}(0)}\right)\right|_{q=0}$$
where $\overrightarrow{\mathbf{F}_{U^+}}(t,s)$ is the operator obtained by quantizing $\displaystyle{p^1_k=\hbar k z^{k} \overrightarrow{\frac{\partial}{\partial q^1_k}}}$ and $\displaystyle{p^2_k=\hbar k \overrightarrow{\frac{\partial}{\partial q^2_k}}}$.\\
Explicitly:
\begin{equation*}
\begin{split}
\mathcal{F}_E \left(\sum t_i \tilde{\Delta}_i+ \sum s_j \tilde{\Theta}_j\right)= & \frac{1}{\hbar}\left[\frac{t_0^2s_2}{2}+t_0t_1s_1- \frac{\hbar s_2}{24}\right]+ \\ + & \mathrm{log}\left[ \mathrm{e}^{\frac{1}{\hbar} \sum \frac{1}{k}\mathrm{e}^{ks_2}\hbar^2 k^2 z^k \overrightarrow{\frac{\partial}{\partial q^1_k}}\overrightarrow{\frac{\partial}{\partial q^2_k}}} \ \mathrm{e}^{\frac{1}{\hbar}\sum \frac{1}{k}q^1_k q^2_k}\right]_{q=0}
\end{split}
\end{equation*}
The differential operator in the second term of the right-hand side can be dealt with by remembering that $\displaystyle{\mathrm{e}^{s \partial_x \partial_y} \mathrm{e}^{xy}=\frac{\mathrm{e}^{\frac{xy}{1-s}}}{1-s}}$, whence (up to an additive constant which is irrelevant for the Gromov-Witten potential)
$$\mathcal{F}_E \left(\sum t_i \tilde{\Delta}_i+ \sum s_j \tilde{\Theta}_j\right)=\frac{1}{\hbar}\left[\frac{t_0^2s_2}{2}+t_0t_1s_1- \hbar \log(\eta(z\mathrm{e}^{s_2}))\right]$$
where $\eta(q)=q^{1/24}\prod_k(1-q^k)$ is the Dedekind eta function.  This agrees with the results of \cite{BCOV}.

\section{Descendent potential of target curves and quantum integrable systems}

In the (quite restricted) family of explicitly computed examples (basically the symplectizations of the low-dimensional contact spheres $S^{2n+1}$), Symplectic Field Theory Hamiltonians show the remarkable property of being integrable, roughly meaning they are part of an infinite number of independent elements in the graded Weyl algebra of $p$'s and $q$'s, forming a graded commutative subalgebra. Indeed, the two simplest cases, $S^1$ and $S^3$, correspond to the (quantum) dispersionless versions of the celebrated KdV and Toda hierarchies of integrable PDEs (see e.g. \cite{DZ}).

Actually this is a general fact (\cite{E}), true for any contact manifold $V$, and the whole commutative algebra of first integrals of the given $\mathbf{H}_{V}$ carries the topological information concerning the so called descendants. More precisely we can extend the definition of the functions $\mathbf{H}$ and $\mathbf{F}$, associated as above to symplectic cobordisms, as generating functions for correlators of the type
\begin{equation*}
\begin{split}
&\langle \tau_{p_1}(\Theta_{i_1}),\ldots,\tau_{p_r}(\Theta_{i_r});\gamma^-_1,\ldots,\gamma^-_{s^-};\gamma^+_1,\ldots,\gamma^+_{s^+} \rangle^d_g =\\
&\int_{\mathcal{M}^A_{g,r,s^-,s^+}} \mathrm{ev}^*(\Theta_{i_1}\otimes\ldots\otimes\Theta_{i_r}\otimes\gamma^-_1\otimes\ldots\otimes\gamma^-_{s^-}\otimes \gamma^+_1\otimes\ldots\otimes\gamma^+_{s^+})\wedge c_1(\psi_1)^{p_1}\wedge\ldots\wedge c_1(\psi_r)^{p_r}
\end{split}
\end{equation*}
where now $c_1(\psi_i)$ is the first Chern class of the tautological line bundle over $\mathcal{M}^A_{g,r,s^-,s^+}$ with fibre at $C \in \mathcal{M}^A_{g,r,s^-,s^+}$ equal to the fibre of $\mathrm{T}^*C$ at the $i$-th marked point. We will call the generating functions of these extended correlators $\mathfrak{H}$ and $\mathfrak{F}$ (the first for the cylindrical case, the second for the general cobordism) and they will depend on the variables $(p,q,t_{i_j,p_k})$ with $t_{i_j,0}=t_{i_j}$ as before.

Theorems \ref{Schrodinger} and \ref{composition} extend then to the case of descendants, i.e. they hold true for the extended functions $\mathfrak{H}$ and $\mathfrak{F}$. Even more fundamentally, the structure equation $[\mathbf{H},\mathbf{H}]=0$ extends to $[\mathfrak{H},\mathfrak{H}]=0$ ($\mathbf{H}$ and $\mathfrak{H}$ are here evaluated on closed forms, so the terms $\di \mathbf{H}$ and $\di \mathfrak{H}$ drop), ensuring that the Hamiltonians appearing in Theorem \ref{Schrodinger} actually commute, forming a quantum integrable system.

\vspace{1cm}

As an example we determine the integrable system involved in the computation of the Gromov-Witten potential of target curves with descendants of $\omega\in H^2(\Sigma_g)$ (in what follows we put even $t_0=0$ for simplicity). We already know that the relevant cylindrical cobordism is just the cylinder $V=S^1\times \R$. Call, as above, $t_0$ and $s_2$ the components along $1$ and $\di \phi$ of a form in $H^\bullet(V)$, then its potential without descendants $\mathbf{H}_0$ is given by (\ref{cylinder}). With a little dimension counting over the moduli space $\mathcal{M}^A_{g,r,s^-,s^+}$ one can also compute the first descendant potential with respect to $\di \phi$ (i.e. the part of $\mathfrak{H}$ depending just on $t_0$ and linearly on $s_{2,1}$) that, with the notation of (\ref{fourier}), takes the form $$s_{2,1}\,\mathbf{H}_1=s_{2,1}\,\frac{1}{\hbar}:\int_{S^1}\left(\frac{u(x)^3}{6}-\hbar \frac{u(x)}{24}\right) \di x:$$
This (quantum) Hamiltonian, corresponding to a quantum dispersionless KdV system, is sufficient to determine uniquely the whole integrable system (notice that it was not so for $\mathbf{H}_0$, since any local Hamiltonian commutes with it), i.e. the algebra of commuting Hamiltonians of Theorem \ref{Schrodinger}. The explicit determination of this commutative algebra has recently been performed in an elegant context of fermionic calculus by Pogrebkov in \cite{P}. Nonetheless the topological information about the target cobordism is encoded also in an explicit choice of a basis of such (super-)commutative algebra.

Such a basis is given (\cite{E}) by
\begin{equation}\label{qdkdvhamiltonians}
\mathbf{H}_{n-2}=\frac{1}{\hbar}\frac{1}{n!}\left.\frac{\partial^n \mathcal{H}(z)}{\partial z^n}\right|_{z=0}
\end{equation}
in terms of the generating function
\begin{equation}\label{qdkdvgenfun}
\mathcal{H}(z)=\frac{1}{\mathcal{S}(\sqrt{\hbar}z)} :\int_{S^1}\mathrm{e}^{z \mathcal{S}(\mathrm{i}\sqrt{\hbar}z\partial_x)\,u(x)}\di x:\qquad\text{with}\qquad\mathcal{S}(t)=\frac{\sinh \frac{t}{2}}{\frac{t}{2}}
\end{equation}
Here $:\cdot:$ means that $q$ and $p$ variables are to be normal ordered ($q$'s are to be put on the left) as if, inside the colon symbols, $q$'s and $p$'s all commuted. As we show in the next section, this generating function can be obtained from the work of Okounkov and Pandharipande, directly by expressing the operator
\begin{equation}\label{diagonal}
\mathcal{E}_0(z)=\sum_{k\in\mathbb{Z}+\frac{1}{2}}\frac{\mathrm{e}^{\sqrt{\hbar}zk}}{\hbar}\, :\psi_k\psi^*_k:+\frac{1}{\mathrm{e}^{\sqrt{\hbar}z/2}-\mathrm{e}^{-\sqrt{\hbar}z/2}}
\end{equation}
of \cite{OP1}, acting on the fermionic Fock space (see below) and relevant for the generating function of $1$-point invariants of $\mathbb{P}^1$ relative to $0$ and $\infty$, in terms of vertex operators via the ``boson-fermion" correspondence (see e.g. \cite{MJD},\cite{K}).

The function $\mathcal{S}(t)$, appearing repeatedly also in the work of Okounkov and Pandharipande, emerges in \cite{P} directly from an integrable systems context (although there a different basis of Hamiltonians for the commutative algebra of symmetries is chosen) and is involved in the topological picture via the Schr\"odinger equation of Theorem \ref{Schrodinger}. For instance, let's write this equation for the simplest case of the descendant SFT-potential for the cap:
\begin{equation}\label{schrodcap}
\mathfrak{F}_\text{cap}(s_{2,k})=\text{log}\left(\mathrm{e}^{\mathfrak{F}_\text{cap}(0)} \prod_{n=0}^{\infty}\mathrm{e}^{s_{2,n}\overleftarrow{\mathbf{H}_{n}}}\right)
\end{equation}
In the next sections, after recalling its fundamental ideas, we make use of fermionic calculus to obtain the fermionic expression of the generating function (\ref{qdkdvgenfun}) and compute explicitly $\mathfrak{F}_\text{cap}(s_{2,k})$.

\section{Fermionic calculus}

Following \cite{MJD}, notice that the Weyl algebra of differential operators of the form $$\sum c_{\alpha_1\alpha_2\ldots\beta_1\beta_2\ldots}p_1^{\alpha_1}p_2^{\alpha_2}\ldots\left(\frac{\partial}{\partial p_1}\right)\left(\frac{\partial}{\partial p_2}\right)\ldots,$$ involved in the Symplectic Field Theory of target curves, is a representation of the algebra $\mathcal{B}$ generated by the \emph{bosons} $na_n:=q_n$ and $a_n^*:=p_n$ (with the usual commutation relations) on the \emph{bosonic Fock space} $\mathbb{C}[[p]]:=\mathbb{C}[[p_1,p_2,\ldots]]=\mathcal{B}\cdot 1$. Notice that the element $1\in\mathbb{C}[[p]]$ (called the \emph{bosonic vacuum state}) is annihilated by $a_n$, for every $n$. In particular $\mathbb{C}[[p]]$ has the basis $$\{a_{m_1}\ldots a_{m_r}\cdot 1\;|\;0<m_1<\ldots<m_r\}.$$

Next we introduce another algebra, the Clifford algebra $\mathcal{A}$ generated by the \emph{fermions} $\{\psi_n,\psi^*_n\}_{m,n\in\mathbb{Z}+1/2}$ with the anti-commutation relations $$\{\psi_m,\psi_n\}=0,\qquad\{\psi^*_m,\psi^*_n\}=0,\qquad\{\psi_m,\psi^*_n\}=\hbar\delta_{m,n}.$$ As for the bosonic case, we represent $\mathcal{A}$ on the \emph{fermionic Fock space} $\mathcal{F}:=\mathcal{A}\cdot|\text{vac}\rangle$, where the \emph{fermionic vacuum state} $|\text{vac}\rangle$ is annihilated by $\psi_k$ with $k<0$ and $\psi^*_l$ with $l>0$. In what follows we will also write $\langle\text{vac}|\cdot\mathcal{A}$ to denote the dual Fock space $\mathcal{F}^*$ (and similarly for its elements).

Besides the bosonic normal ordering of $p$ and $q$ operators, which we already defined, we introduce the following fermionic normal ordering:
\begin{equation*}
:\psi_m \psi^*_n:=\begin{cases}\psi_m\psi^*_n & \text{if}\ m>0\ \text{or if}\ n>0 \\ -\psi^*_n\psi_m & \text{if}\ m<0\ \text{or if}\ n<0 \end{cases}
\end{equation*}

The space $\mathcal{F}$ (and its dual) is graded according to the so called \emph{charge} $l$ as $\displaystyle{\mathcal{F}=\bigoplus_{l\in\mathbb{Z}} \mathcal{F}_l}$ where $\mathcal{F}_l$ is the span of $$\{\psi_{m_r}\ldots\psi_{m_1}\psi^*_{n_1}\ldots\psi^*_{n_s}|\text{vac}\rangle\;| \;m_1<\ldots<m_r,\,n_1<\ldots<n_s<0,\,r-s=l\}.$$

\vspace{1cm}

The following is an important and well known result (proven e.g. in \cite{MJD}) usually called the boson-fermion correspondence.

\begin{thm}
There is an isomorphism $\Phi:\mathcal{F}_0\to\mathbb{C}[[p]]$ of the form $$\Phi(|u\rangle):=\langle\mathrm{vac}|\mathrm{e}^{K(p)}|u\rangle.$$
where $K(p):=\frac{1}{\hbar}\sum_{n=1}^\infty \frac{p_n}{n}\sum_{j\in\mathbb{Z}+1/2}:\psi_j\psi^*_{j+n}:$\\
At the level of operators the isomorphism is described by 
$$\Phi^{-1}q_n\Phi=\sum_{j\in\mathbb{Z}+1/2}:\psi_j\psi^*_{j+n}:$$
$$\Phi^{-1}p_n\Phi=\sum_{j\in\mathbb{Z}+1/2}:\psi_j\psi^*_{j-n}:$$
\end{thm}

\begin{example}\label{fermionex}
Using the formula $\mathrm{e}^A B \mathrm{e}^{-A}=\mathrm{e}^{\mathrm{ad}_A} B$ for the adjoint action of a Lie algebra on itself we get
$$\mathrm{e}^{K(p)}\psi_n\mathrm{e}^{-K(p)}=\psi_n+\frac{p_1}{\hbar}\psi_{n-1}+ \left(\frac{p_2}{\hbar}+\frac{p_1^2}{2\hbar^2}\right)\psi_{n-2}+\ldots$$
$$\mathrm{e}^{K(p)}\psi^*_n\mathrm{e}^{-K(p)}=\psi^*_n-\frac{p_1}{\hbar}\psi^*_{n+1}+ \left(-\frac{p_2}{\hbar}+\frac{p_1^2}{2\hbar^2}\right)\psi^*_{n+2}+\ldots$$
and, using these, we can compute for instance 
\begin{equation*}
\begin{split}
\Phi(\psi_{3/2}\psi^*_{-1/2}|\mathrm{vac}\rangle)=& \langle\mathrm{vac}|\mathrm{e}^{K(p)}\psi_{3/2}\psi^*_{-1/2}|\mathrm{vac}\rangle\\
=&\langle\mathrm{vac}|\left(\psi_{3/2}+\frac{p_1}{\hbar}\psi_{1/2}+ \left(\frac{p_2}{\hbar}+\frac{p_1^2}{2\hbar^2}\right)\psi_{-1/2} +\ldots\right)\left(\psi^*_{-1/2}+\ldots\right)|\mathrm{vac}\rangle\\
=&\frac{p_2}{2}+\frac{p_1^2}{2\hbar}
\end{split}
\end{equation*}
and similarly $$\Phi(\psi_{1/2}\psi^*_{-3/2}|\mathrm{vac}\rangle)=\frac{p_2}{2}-\frac{p_1^2}{2\hbar}$$
\end{example}

Now consider the following \emph{vertex operator} acting (on the right) on $\mathbb{C}[[p]]$
$$X(z_1,z_2):=\mathrm{e}^{\sum_{j=1}^\infty\sqrt{\hbar}(z_1^j-z_2^j)\overleftarrow{\frac{\partial}{\partial p_j}}}\, \mathrm{e}^{-\sum_{j=1}^\infty\frac{1}{\sqrt{\hbar}j}(z_1^{-j}-z_2^{-j})p_j}$$
with $z_1,z_2 \in \mathbb{C}$. Then the above theorem gives as a corollary
\begin{equation}\label{vertex}
\frac{1}{\hbar}:\psi(z_1)\psi^*(z_2):\ \xrightarrow{\Phi}\ \frac{1}{z_1-z_2}\left(X(z_1,z_2)-1\right)
\end{equation}
where $\psi(z)=\sum_{n\in\mathbb{Z}+1/2}\psi_n z^{n-1/2}$ and $\psi^*(z)=\sum_{n\in\mathbb{Z}+1/2}\psi^*_n z^{-n-1/2}$

\vspace{1cm}

Getting back to the generating function (\ref{qdkdvgenfun}) we make the following important observation: 
$$\overleftarrow{\mathcal{H}(z)}=\frac{1}{\mathcal{S}(\sqrt{\hbar}z)}\int_{S^1} X\left(\mathrm{e}^{\mathrm{i}\left(x-\mathrm{i}\frac{\sqrt{\hbar}z}{2}\right)}, \mathrm{e}^{\mathrm{i}\left(x+\mathrm{i}\frac{\sqrt{\hbar}z}{2}\right)}\right)\mathrm{d}x$$
whence, by applying (\ref{vertex}) and computing the residue, we get the fermionic expression for the generating function of the quantum dispersionless KdV Hamiltonians
\begin{thm}
\begin{equation}\label{fermionicgenfun}
\overleftarrow{\mathcal{H}(z)}\ \xrightarrow{\Phi}\ \sqrt{\hbar} z \mathcal{E}_0(z) =\sum_{k\in\mathbb{Z}+\frac{1}{2}}\frac{\sqrt{\hbar} z\mathrm{e}^{\sqrt{\hbar}zk}}{\hbar}\, :\psi_k\psi^*_k:+\frac{\sqrt{\hbar} z}{\mathrm{e}^{\sqrt{\hbar}z/2}-\mathrm{e}^{-\sqrt{\hbar}z/2}}
\end{equation}
\end{thm}
This way the equivalence of our generating function of $1$-point SFT-invariants of the tube and the one for $1$-point GW-invariants relative to $0$ and $\infty$ (derived by Okounkov and Pandharipande in \cite{OP1} using the Gromov-Witten/Hurwitz correspondence), is proven. From (\ref{fermionicgenfun}) we get
\begin{equation}\label{fermionicham}
\overleftarrow{\mathbf{H}_n} \ \xrightarrow{\Phi}\ \frac{\hbar^{n/2}}{(n+1)!}\left[\sum_{k\in\mathbb{Z}+1/2}\frac{k^{n+1}}{\hbar}:\psi_k\psi^*_k: +(1-2^{-(n+1)})\zeta(-(n+1))\right]
\end{equation}
where $\zeta$ is the Riemann zeta function. The peculiar diagonal form of the right-hand side ensures the required commutativity of the Hamiltonians. Now we can use this fermionic formalism to compute the invariants of $\mathbb{P}^1$ (or, in principle, of any compact Riemann surface) via Theorems \ref{Schrodinger} and \ref{composition}.

\section{Descendant potential of $\mathbb{P}^1$}

First of all we rewrite equation (\ref{schrodcap}) for the descendant SFT-potential of the cap in the language of fermionic calculus of the previous section. We get
$$\mathrm{e}^{\mathfrak{F}_\mathrm{cap}(s_{2,k})}=\mathrm{e}^{\sum_n s_{2,n}\frac{\hbar^{n/2}}{(n+1)!}\left[\sum_{k\in\mathbb{Z}+1/2}\frac{k^{n+1}}{\hbar}:\psi_k\psi^*_k: +(1-2^{-(n+1)})\zeta(-(n+1))\right]}\ \mathrm{e}^{\sum_{j\in\mathbb{Z}+1/2}:\psi_j\psi^*_{j-1}:}\ |\mathrm{vac}\rangle$$
One is now to compute the right hand side and then use Theorem \ref{composition} to get the descendant potential of $\mathbb{P}^1$. If we restrict to the case of maps with fixed degree the computations can be carried out very explicitly. Remember that both the SFT and Gromov-Witten potentials are power series in $z$ with coefficients the potentials for maps with fixed degree. We write
$$\mathcal{F}_{\mathbb{P}^1}(s_{2,k},z)=\sum_i \mathcal{F}_{\mathbb{P}^1,i}(s_{2,k}) z^i$$
$$\mathfrak{F}_{\mathrm{cap}}(s_{2,k},p\mapsto zp)=\sum_i \mathfrak{F}_{\mathrm{cap},i}(s_{2,k},p) z^i$$
Then Theorems \ref{Schrodinger} and \ref{composition} give:\\
\\
\textsc{degree} $0$: $$\mathrm{e}^{\mathcal{F}_{\mathbb{P}^1,0}(s_{2,k})}=\mathrm{e}^{\overrightarrow{\mathfrak{F}}_{\mathrm{cap},0}(s_{2,k},p)}\ \mathrm{e}^{\mathfrak{F}_{\mathrm{cap}}(0,q)}$$
$$\mathrm{e}^{\mathfrak{F}_{\mathrm{cap},0}(s_{2,k},p)}=1\ \mathrm{e}^{\sum_{n}s_{2,n}\overleftarrow{\mathbf{H}}_n}$$
and using (\ref{fermionicham}) we immediately get (see \cite{FP},\cite{Pa})
$$\mathfrak{F}_{\mathrm{cap},0}(s_{2,k},p)=\sum_{n}s_{2,n}\frac{\hbar^n}{(n+1)!}(1-2^{-(n+1)})\zeta(-(n+1))$$
$$\mathcal{F}_{\mathbb{P}^1,0}(s_{2,k})=\sum_{n}s_{2,n}\frac{\hbar^n}{(n+1)!}(1-2^{-(n+1)})\zeta(-(n+1))$$
\\
\textsc{degree} $1$: $$\mathcal{F}_{\mathbb{P}^1,1}(s_{2,k})=\overrightarrow{\mathfrak{F}}_{\mathrm{cap},1}(s_{2,k},p)\ \mathrm{e}^{\mathfrak{F}_{\mathrm{cap}}(0,q)}$$
$$\mathfrak{F}_{\mathrm{cap},1}(s_{2,k},p)\ \mathrm{e}^{\mathfrak{F}_{\mathrm{cap},0}(s_{2,k},p)}=\frac{1}{\hbar}p_1\ \mathrm{e}^{\sum_{n}s_{2,n}\overleftarrow{\mathbf{H}}_n}$$
and passing to the fermions the second equation becomes
$$\Phi^{-1}(\mathfrak{F}_{\mathrm{cap},1}(s_{2,k},p))= \frac{1}{\hbar}\mathrm{e}^{\sum_k\left(\sum_n s_{2,n}\frac{\hbar^{(n-2)/2}}{(n+1)!}k^{n+1}\right):\psi_k\psi^*_k:}\ \psi_{1/2}\psi_{-1/2}\ |\mathrm{vac}\rangle$$
Since $$\left[\sum_k\left(\sum_n s_{2,n}\frac{\hbar^{(n-2)/2}}{(n+1)!}k^{n+1}\right):\psi_k\psi^*_k:\ ,\ \psi_{1/2}\psi_{-1/2}\right]=\sum_n s_{2,2n}\frac{\hbar^n}{2^{2n}(2n+1)!}\ \psi_{1/2}\psi_{-1/2},$$ we get
\begin{thm}
$$\mathfrak{F}_{\mathrm{cap},1}(s_{2,k},p)=\frac{1}{\hbar}\mathrm{e}^{\sum_n s_{2,2n}\frac{\hbar^n}{2^{2n}(2n+1)!}}\ p_1$$
$$\mathcal{F}_{\mathbb{P}^1,1}(s_{2,k})=\frac{1}{\hbar}\mathrm{e}^{\sum_n s_{2,2n}\frac{\hbar^n}{2^{2n}(2n+1)!}}$$
\end{thm}
The second equation coincides with the result found by Pandharipande in \cite{Pa} using Toda equations.
\\
\textsc{degree} $2$:
$$\left(\mathcal{F}_{\mathbb{P}^1,1}(s_{2,k})\right)^2+2\mathcal{F}_{\mathbb{P}^1,2}(s_{2,k})= \left[\left(\overrightarrow{\mathfrak{F}}_{\mathrm{cap},1}(s_{2,k},p)\right)^2+ 2\overrightarrow{\mathfrak{F}}_{\mathrm{cap},2}(s_{2,k},p)\right]\ \mathrm{e}^{\mathfrak{F}_{\mathrm{cap}}(0,q)}$$
$$\left[\left(\mathfrak{F}_{\mathrm{cap},1}(s_{2,k},p)\right)^2+ 2\mathfrak{F}_{\mathrm{cap},2}(s_{2,k},p)\right]\ \mathrm{e}^{\mathfrak{F}_{\mathrm{cap},0}(s_{2,k},p)}=\frac{1}{\hbar^2}p_1^2\ \mathrm{e}^{\sum_{n}s_{2,n}\overleftarrow{\mathbf{H}}_n}$$
and passing to fermions the second equation becomes (see Example \ref{fermionex})
\begin{equation*}
\begin{split}
&\Phi^{-1}(\left(\mathfrak{F}_{\mathrm{cap},1}(s_{2,k},p)\right)^2+ 2\mathfrak{F}_{\mathrm{cap},2}(s_{2,k},p))=\\
&\frac{1}{\hbar}\mathrm{e}^{\sum_k\left(\sum_n s_{2,n}\frac{\hbar^{(n-2)/2}}{(n+1)!}k^{n+1}\right):\psi_k\psi^*_k:}\ (\psi_{3/2}\psi^*_{-1/2}-\psi_{1/2}\psi^*_{-3/2})\ |\mathrm{vac}\rangle
\end{split}
\end{equation*}
Since
\begin{equation*}
\begin{split}
&\left[\sum_k\left(\sum_n s_{2,n}\frac{\hbar^{(n-2)/2}}{(n+1)!}k^{n+1}\right):\psi_k\psi^*_k:\ ,\  (\psi_{3/2}\psi^*_{-1/2}-\psi_{1/2}\psi^*_{-3/2})\right]=\\
&\sum_n s_{2,n}\frac{\hbar^n/2}{(n+1)!}\left(\left(\frac{3}{2}\right)^{n+1}-\frac{1}{(-2)^{n+1}}\right) \psi_{3/2}\psi^*_{-1/2}\\
&-\sum_n s_{2,n}\frac{\hbar^n/2}{(n+1)!}\left(\frac{1}{2^{n+1}} -\left(\frac{3}{2}\right)^{n+1}\right) \psi_{1/2}\psi^*_{-3/2},
\end{split}
\end{equation*}
we get
\begin{equation*}
\begin{split}
&\left(\mathfrak{F}_{\mathrm{cap},1}(s_{2,k},p)\right)^2+ 2\mathfrak{F}_{\mathrm{cap},2}(s_{2,k},p)=\\
&\frac{1}{\hbar}\mathrm{e}^{\sum_n s_{2,n}\frac{\hbar^n/2}{(n+1)!}\left(\left(\frac{3}{2}\right)^{n+1}- \frac{1}{(-2)^{n+1}}\right)}\left(p_2+\frac{p_1^2}{2\hbar}\right)
-\frac{1}{\hbar}\mathrm{e}^{\sum_n s_{2,n}\frac{\hbar^n/2}{(n+1)!}\left(\frac{1}{2^{n+1}} -\left(\frac{3}{2}\right)^{n+1}\right)}\left(p_2-\frac{p_1^2}{2\hbar}\right)
\end{split}
\end{equation*}
and finally
\begin{thm}
\begin{equation*}
\begin{split}
&\mathcal{F}_{\mathbb{P}^1,2}(s_{2,k})=\\
&\frac{1}{2\hbar^2}\left[\frac{1}{2}\mathrm{e}^{\sum_n s_{2,n}\frac{\hbar^n/2}{(n+1)!}\left(\left(\frac{3}{2}\right)^{n+1}- \frac{1}{(-2)^{n+1}}\right)}
-\frac{1}{2}\mathrm{e}^{\sum_n s_{2,n}\frac{\hbar^n/2}{(n+1)!}\left(\frac{1}{2^{n+1}} -\left(\frac{3}{2}\right)^{n+1}\right)}
+\mathrm{e}^{\sum_n s_{2,2n}\frac{\hbar^n}{2^{2n-1}(2n+1)!}}\right]
\end{split}
\end{equation*}
\end{thm}

This process can be carried on to even higher degree, with more and more struggle, but it appears, already at degree $2$, to be fairly more efficient than the method of Virasoro constraints of \cite{EHX}, used in \cite{OP2}.

\end{document}